 \documentclass[jaa]{degruyter-journal-a}      
\usepackage{mathrsfs}

\theoremstyle{plain} 
\newtheorem{cor}{Corollary}

\newtheorem{thm}{Theorem}

\theoremstyle{definition}
\newtheorem{defn}{Definition}

\theoremstyle{remark}

\numberwithin{equation}{section}

\DeclareMathOperator{\divgce}{Div}
\DeclareMathOperator{\pd}{\partial}

\DeclareFontFamily{OT1}{pzc}{}
\DeclareFontShape{OT1}{pzc}{m}{it}{<-> s * [1.10] pzcmi7t}{}
\DeclareMathAlphabet{\mathpzc}{OT1}{pzc}{m}{it}

\DeclareMathSymbol{\R}{\mathalpha}{AMSb}{"52}
\DeclareMathSymbol{\C}{\mathalpha}{AMSb}{"43}
\newcommand{\mbb}[1]{\mathbb{#1}}
\newcommand{\F}{\mbb{F}}

\newcommand{\mscl}{\mathscr{L}}
\newcommand{\mscs}{\mathscr{S}}

\newcommand{\bv}{\mathbf{v}}

\newcommand{\fq}{\mathfrak{q}}

\newcommand{\mcs}{\mathcal{S}}
\newcommand{\mcl}{\mathcal{L}}

\newcommand{\set}[1]{\left\{#1\right\}}

\title{Some variational principles  associated with ODEs of maximal symmetry. Part 2: The general case}
\headlinetitle{Variational principles  associated with ODEs of maximal symmetry}

\lastnameone{Ndogmo}
\firstnameone{Jean-Claude}
\nameshortone{J.\,C.~Ndogmo}
\addressone{Department of Mathematics and Applied Mathematics,
University of Venda,
P/B X5050, Thohoyandou,
Limpopo  0950}
\countryone{South Africa}
\emailone{ndogmoj@yahoo.com}

\abstract{Variational and divergence symmetries are studied  in this paper for  the whole class of linear and nonlinear equations of maximal symmetry, and the associated first integrals are given in explicit form.  All the main results obtained  are formulated as theorems or conjectures for equations of a general order. A discussion of the existence of variational symmetries with respect to a different Lagrangian, which turns out to be the most common and most readily available one, is also carried out. This leads to significantly different results when compared with the former case of the transformed Lagrangian. The latter analysis also gives rise to more general results concerning the variational symmetry algebra of any linear or nonlinear equations.\par}



\researchsupported{This research is supported by  research grants from the University of Venda and from the NRF of South Africa.}

\begin{document}

\section{Introduction}
\label{s:intro}

An important class of nonlinear ordinary differential equations ({\small\sc ode}s) consists of those equations of maximal symmetry, which means having a Lie point symmetry algebra of maximal dimension. Although by a result of Lie \cite{LieCanonic, olv94} any such nonlinear {\small \sc ode} can be mapped by a point transformation to the trivial equation, there is however no easy means to identify any of them and finding the point transformation mapping the nonlinear equation to the linear counterpart involves solving a system of partial differential equations which in certain cases might turn out to be more difficult than solving the nonlinear equation itself, if not impossible. On the other hand, many important properties of linear equations of maximal symmetry can be more easily found and could possibly be transferred to their nonlinear counterparts.\par

Due to some relatively recent results on linear equations of  maximal symmetry algebra \cite{KM, ML, JF, J16} it should be possible by now to describe in their most general form the variational principles associated with this class of equations and with the whole class of nonlinear  {\small \sc ode}s of maximal symmetry, and to obtain in particular the precise conditions for the existence of variational symmetries, and the related first integrals thereof. \par

Pursuing a work undertaken by Lie \cite{LieCanonic}, Krause and Michel \cite{KM} have indeed given a specific proof of the fact that a linear ordinary differential equation ({\small \sc lode}) has maximal symmetry if and only if it can be reduced by a point transformation to the form $y^{(n)}=0,$ which will henceforth be referred to as the canonical form.

It is well known that a {\small \sc lode} can always be reduced by a point transformation to the normal form
\begin{equation} \label{nor1}
y^{(n)} + A_n^2 y^{(n-2)}+ \dots + A_n^{n-1} y^{(1)} + A_n^{n} y =0,
\end{equation}
in which the coefficient of the term of second highest order has vanished. A specific property satisfied by {\small \sc lode}s of maximal symmetry in the form \eqref{nor1} is that $n$ linearly independent solutions $s_k$ can be obtained in the form
\begin{equation}\label{sk}
s_k=  u^{n-k-1} v^k, \qquad 0 \leq k \leq n-1, \qquad n \geq 2,
\end{equation}
for some functions $u$ and $v$ of $x.$ Moreover, when such equations are in normal form \eqref{nor1}, the Wronskian  $\mathpzc{w}(u,v)= \det {\tiny \begin{pmatrix} u & v\\ u'  & v' \end{pmatrix} }$  is a nonzero constant and will be normalized to one. It also follows from \eqref{sk} that $u$ and $v$ are two linearly independent solutions of the second-order source equation
\begin{equation} \label{srce}
y_{xx} + \fq(x)y=0,\qquad \fq= A_2^2.
\end{equation}
Using in particular this fact, the symmetry generators of equations of the form \eqref{nor1} of maximal symmetry were all obtained in \cite{KM} in terms of two linearly independent solutions $u$ and $v$  of \eqref{srce}. In the sequel, when we refer to \eqref{nor1}, we shall assume that it is  the most general form of a {\small \sc lode} of maximal symmetry, unless otherwise stated.\par

In this paper, by exploiting the expression of the symmetry generators for equations of the form \eqref{nor1} of maximal symmetry, we investigate for these equations and for the whole class of nonlinear  {\small \sc ode}s of maximal symmetry, the existence, and carry out the determination of  variational and  divergence symmetry algebras, and first integrals in the sense given in \cite{olv1}. This study is an application of results obtained for equations in canonical form in Part 1, see \cite{part1},  of this series of two papers. \par

In the case of linear  {\small \sc ode}s, given that analytical expressions for the general coefficient   of {\small \sc lode}s of maximal symmetry is not yet available for equations in the general form \eqref{nor1}, some of the results in the study of \eqref{nor1} will often be limited to equations of low order not exceeding eight for which  explicit expressions are available \cite{ML, JF}. Nevertheless,  all the main results obtained are extended  to general theorems  or to natural conjectures for linear and nonlinear {\small \sc ode}s of a general order. \par

A discussion of the existence of variational symmetries with respect to a different Lagrangian, which turns out to be the most common and most readily available one, is also carried out. This leads to significantly different results when compared with the former case of the transformed Lagrangian. The latter analysis also gives rise to more general results concerning the variational symmetry algebra of any linear or nonlinear equations.\par

This short paper is to be read in conjunction with Part 1, see \cite{part1},  of this series two papers. In particular, some of the equation numbers referenced in this paper refer to those introduced in \cite{part1}. The field $\F$  of scalars will be assumed to be the field $\R$ of real numbers, although $\F$ could equally be taken to be $\C$.

\section{Preliminary results}
\label{s:basic}

We collect in this section some of the preliminary results needed for our discussion of variational symmetries and first integrals. The relevant facts concerning variational symmetries appear in \cite[Section 2]{part1} of this series of two papers, and should be read in conjunction with the actual paper, as already stated. In particular, the notations related to variational symmetries to be used are also those defined in \cite[Section 2]{part1}.\par

For  a given Lagrangian $L$  corresponding to an {\sc ode} $\Delta_n[y]\equiv \Delta=0,$ set
\begin{align}
\mathscr{S}(\mathbf{v}) &=  \mathbf{v}^{[n]} (L) + L \divgce \xi,\qquad (\text{ for $n$ even}) \label{plv} \\
\mathscr{D}(\mathbf{v}) &= D_\Delta^* (Q) + D_Q^* (\Delta), \label{adjc}
\end{align}
where $\mathbf{v}= \xi(x,y) \pd_x + \psi(x,y)\pd_y$ has characteristic $Q.$ Then a (Lie point) symmetry $\bv$ of $\Delta=0$ is a variational symmetry if and only if
\begin{subequations}\label{v&d-cdt}
\begin{align}
\mathscr{S}(\mathbf{v})=& 0, \label{vsymeq} \\[-3mm]
\intertext{and a divergence symmetry   if and only if}\vspace{-7mm}
\mathscr{D}(\mathbf{v}) =&0  \label{dsymeq}.\vspace{-5mm}
\end{align}
\end{subequations}

On the other hand,  $n+4$ infinitesimal generators for an $n$th order {\small \sc lode} of the general form \eqref{nor1} of maximal symmetry are given by
\begin{subequations} \label{ig}
\begin{align}
V_k &= s_k \pd_y = u^{n-(k+1)} v^k \pd_y, \qquad 0\leq k \leq n-1 \label{ig1}\\
W_y &= y \pd_y    \label{ig2}\\
F_n &= u^2 \pd_x + (n-1) u u' y \pd_y    \label{ig3}\\
G_n &= 2 u v \pd_x + (n-1) (u v' + u' v) y \pd_y  \label{ig4}\\
H_n &= - v^2 \pd_x - (n-1) v v' y \pd_y,    \label{ig5}
\end{align}
\end{subequations}

where the $s_k$ are as in \eqref{sk} $n$ linearly independent solutions of the $n$th order equation, while $u$ and $v$ are two linearly independent solutions of the second order source equation \eqref{srce}.\par

It should also be mentioned here that the symmetries $V_k$  are often referred to as  solution symmetries and they generate the $n$-dimensional abelian Lie algebra $\mathcal{A}_n,$ while $W_y$ is often called the homogeneity symmetry. On the other hand, the Lie algebra $\mathfrak{s}$ generated by $F_n, G_n$ and $H_n$ is isomorphic to $\mathfrak{sl}_2 \equiv \mathfrak{sl}(2,\F)$ and the full symmetry algebra $\mathfrak{g}_n$ of \eqref{nor1} is the semi-direct sum $\mathfrak{g}_n= (\mathcal{A}_n \dotplus \F W_y) \ltimes \mathfrak{s},$ where the operator $\ltimes$ denotes a semi-direct sum of Lie algebras. We shall often denote by $\langle X_1,\dots,X_r \rangle$ a Lie algebra generated by the  vectors $X_1,\dots X_r.$ Thus for instance $\mathcal{A}_n= \langle V_0, \dots, V_{n-1} \rangle. $


\section{Equations in the most general form }\label{s:gform}

We now proceed to the direct investigation of variational and divergence symmetries for {\small \sc lode}s  of maximal symmetry in their most general form \eqref{nor1}. Just as equivalent equations (under point transformations) have isomorphic Lie point symmetry algebras, they also have isomorphic divergence  symmetry algebras.  In the case of lagrangian equations, not only they are not invariant under point transformations, but also two Lagrangian equations equivalent by point transformation need not have  isomorphic variational symmetry algebras, except if  the Lagrangian considered for the transformed equation is the transformed counterpart of the Lagrangian of the original equation under the change of coordinates. These results can be formulated more formally as follows. We first give a more formal definition of a Lagrangian equation.\par
\begin{defn}
A differential equation $\Delta=0$ is called Lagrangian if it is equivalent to an Euler--Lagrange equation $E (L)=0,$ with Lagrangian function $L.$ Here, two equations are said to be equivalent if they have the same solution set.
\end{defn}

Denote by $Z$ and $X$ the space $\R^p$ of independent variables coordinatized by  $z= (z^1, \dots, z^p)$ and $x= (x^1,\dots, x^p),$ respectively. Similarly, denote by $W$ and $Y$ the space $\R^q$ of dependent variables coordinatized by $y= (y^1, \dots, y^q)$ and $w= (w^1, \dots, w^q),$ respectively.
\begin{thm}\label{th:transf-vsym}
With the above notation, let
\begin{equation} \label{pt}
z= \zeta(x, y),\qquad \text{ and }\qquad  w= \phi(x, y)
\end{equation}
be an invertible change of variables defined from an open subset $\Omega \subset X \times Y$ onto an open subset $\tilde{\Omega} \subset Z\times W,$ and mapping the (linear or nonlinear) differential equation $\tilde{\Delta} \equiv \tilde{\Delta} (z, w^{(n)})=0$ to $\Delta \equiv \Delta(x, y^{(n)})=0.$
Let $\tilde{\bv}$ be a vector field defined on $\tilde{\Omega}$ and  $\bv$  its transformed version under \eqref{pt}.
\begin{enumerate}
\item[(a)] $\tilde{\bv}$ is a divergence symmetry for $\tilde{\Delta}=0$ if and only if $\bv$ is a divergence symmetry for $\Delta=0.$

\item[(b)] Suppose that $\tilde{\Delta}=0$ is a Lagrangian equation with Lagrangian $\tilde{L},$ and that  $\Delta=0$ is also a Lagrangian equation. Then the Lagrangian $L$ for $\Delta=0$  is the transformed version of  $\tilde{L}$ under to the change of coordinates \eqref{pt}. In this case, $\tilde{\bv}$ is a variational symmetry for $\tilde{\Delta}=0$ with corresponding Lagrangian $\tilde{L}$ if and only if $\bv$ is a variational symmetry of $\Delta=0$ with corresponding Lagrangian $L.$
\end{enumerate}
\end{thm}
The point transformations mapping the trivial equation $w^{(n)}=0$ to the most general form of a linear equation of maximal symmetry  of the form of \cite[(3.1)]{part1} in which it depends on a single arbitrary function $\fq$ is well-known \cite{J16}. These transformations are given by
\begin{equation} \label{ptmaxsym}
z= \int \frac{1}{u^2} dx ,\qquad \text{ and }\qquad  w= \alpha u^{1-n} y
\end{equation}
where $u$ is as usual a nonzero solution of the second-order source equation \eqref{srce}, and $\alpha$ is an arbitrary constant. It therefore follows from   \cite[Conjecture 1]{part1}  and Theorem \ref{th:transf-vsym} that the divergence symmetry algebras and the variational symmetry algebras for linear equations of maximal symmetry of the most general form can be found. Before describing  these subalgebras, we first need to find a precise expression of the transformed Lagrangian of the form \cite[(4.9)]{part1}  under the change of coordinates \eqref{ptmaxsym}. However, contrary to the case of equations in canonical form $y^{(n)}=0$ the general expression of the  transformed version $\mathscr{L}_n$ under \eqref{ptmaxsym} of the Lagrangian \cite[(4.9)]{part1}  is not available in closed form for arbitrary orders $n.$ The corresponding expression of  $\mathscr{L}_n$ for general equations  \cite[(3.1)]{part1}  of even orders $n$ is given for an $m$th-order Lagrangian $L_n(z , w^{(m)})$ by
\begin{equation} \label{lagran-trf}
\mathscr{L}_n (x, y^{(m)}) = L_n(z , w^{(m)})\, z_{\, x}.
\end{equation}
For instance, the expressions for $\mathscr{L}_n$ such that $2\leq n\leq 6$ are given as follows, where $\mathfrak{i}= u'/u.$
\begin{subequations} \label{tfLn}
\begin{align}
\mathscr{L}_2=& -\frac{1}{2} \mathfrak{i}^2 y^2+\mathfrak{i} y y_x-\frac{y_x^2}{2} \label{tfL2}\\
\begin{split}
\mathscr{L}_4=& \frac{9}{2} \left(\fq+2 \mathfrak{i}^2\right)^2 y^2-12 \mathfrak{i} \left(\fq+2 \mathfrak{i}^2\right) y y_x+8 \mathfrak{i}^2 y_x^2\\
&+3 \left(\fq+2 \mathfrak{i}^2\right) y y_{xx}-4 \mathfrak{i} y_x y_{xx}+\frac{y_{xx}^2}{2}
\end{split}  \label{tfL4}\\
\begin{split}
\mathscr{L}_6=& \frac{25}{2} y^2 \left(9 \fq \mathfrak{i}+12 \mathfrak{i}^3-\fq_x\right)^2+5 \left(13 \fq+36 \mathfrak{i}^2\right) y \left(9 \fq \mathfrak{i}+12 \mathfrak{i}^3-\fq_x\right) y_x \\
&-\frac{1}{2} \left(13 \fq+36 \mathfrak{i}^2\right)^2 y_x^2-45 \mathfrak{i} y \left(9 \fq \mathfrak{i}+12 \mathfrak{i}^3-\fq_x\right) y_{xx}\\
&+9 \left(13 \fq \mathfrak{i}+36 \mathfrak{i}^3\right) y_x y_{xx}-\frac{81}{2} \mathfrak{i}^2 y_{xx}^2+y \left(45 \fq \mathfrak{i}+60 \mathfrak{i}^3-5 \fq_x\right)y^{(3)}\\
&+\left(-13 \fq-36 \mathfrak{i}^2\right) y_x y^{(3)}+9 \mathfrak{i} y_{xx} y^{(3)}-\frac{1}{2} {y^{(3)}}^2
\end{split}  \label{tfL6}
\end{align}
\end{subequations}

The result about the variational and the divergence symmetry algebras of the transformed equations can now be stated as follows.
\begin{cor} \label{cor:gen-vsym}
Denote by $\Delta_n\equiv\Delta(x, y^{(n)})=0$ the transformed version of the trivial equation $w^{(n)}=0$ under the change of coordinates \eqref{ptmaxsym}.  Thus $\Delta_n=0$ is the most general form \eqref{nor1} of a linear equation of maximal symmetry  and depends on an arbitrary function $\fq=\fq(x).$ For even orders $n,$ denote by $\mathscr{L}_n$ the transformed version of the Lagrangian $L_n$ of $w^{(n)}=0$ as given by  \cite[(4.9)]{part1} . Denote by $\mathcal{S}_{div}$ and $\mathcal{S}_{var}$ the divergence symmetry algebra and the variational symmetry algebra of $\Delta_n=0.$

\begin{enumerate}
\item[(a)] For $n$ odd, $\mathcal{S}_{div}= \langle V_0, \dots,  V_{n-1}, W_y\rangle.$ That is, $\mathcal{S}_{div}= \mathcal{A}_n \dotplus \F W_y.$
\item[(b)] For $n$ even, $\mathcal{S}_{div} = \mathcal{A}_n \dotplus \mathfrak{s},$ and with respect to the Lagrangian $\mathscr{L}_n$ one has $\mathcal{S}_{var} =  \langle V_0, \dots,  V_{\frac{n-2}{2}}\rangle \dotplus  \langle F_n, G_n\rangle.$
\end{enumerate}
\end{cor}
\begin{proof}
This is an immediate consequence of \cite[Conjecture 1]{part1}  and Theorem \ref{th:transf-vsym}.
\end{proof}

It would be desirable to also give a complete determination of first integrals associated with equations of maximal symmetry and depending on an arbitrary function $\fq,$ as it was done in  \cite[Section 4]{part1}  with $\fq=0$. However, this problem can only be briefly discussed here due to space limitations. \par
The good news about the first integral of the transformed equation is that they can all be obtained from those of the original equation by the mere change of coordinates defined by the related point transformation. However, in the case of the transformed {\small \sc lode}s $\Delta_n[y]=0,$ the point transformation \eqref{ptmaxsym} mapping $y^{(n)}=0$ to them depends on a solution  $u$ of the source equation. Therefore, the transformed first integral will generally also depend on $u.$ This could be a drawback because first integrals are often used for the determination of solutions of differential equations or to study the properties such solutions.\par

With some additional efforts, these solutions can be eliminated from the expressions of the transformed first integrals. We consider for instance the case of the homogeneity symmetry $W_y$ which has an invariant expression for all linear equations. Although the expression for $W_y$ does not depend on $n,$ $\mathscr{F} (W_y)$ does, and thus for the sake of clarity we denote by $\mathscr{F}^n (W_y)$  the first integral associated with the homogeneity symmetry, for each order $n.$ Then by applying the change of variables \eqref{ptmaxsym} to \cite[(4.15)]{part1}  and eliminating any solution in the resulting expression, one obtains the following expressions for $n=3,5, 7.$
\begin{align} \label{fi-wy-gen}
\begin{split}
\mathscr{F}^3 (W_y) &= 2 \fq y^2-\frac{y_x^2}{2}+y y_{xx}  \\
 \mathscr{F}^5 (W_y) &= 10 y \fq_x y_x-10 \fq y_x^2+4 y^2 \left(8 \fq^2+\fq_{x x}\right)\\
&\quad +20 \fq y y_{xx}+\frac{y_{xx}^2}{2}-y_x y^{(3)}+y y^{(4)} \\
\mathscr{F}^7 (W_y) &=  -28 y_x^2 \left(14 \fq^2+\fq_{x x}\right)-28 \fq_x y_x y_{xx}\\
&\quad +y \left(784 \fq^2+84 \fq_{x x}\right) y_{xx}+28 \fq y_{xx}^2\\
&\quad +28 y y_x \left(28 \fq \fq_x+\fq^{(3)}\right)+84 y \fq_x y^{(3)}-56 \fq y_x y^{(3)}-\frac{1}{2} y_{3}^2\\
&\quad +6 y^2 \left(192 \fq^3+33 \fq_x^2+52 \fq \fq_{x x}+\fq^{(4)}\right)+56 \fq y y^{(4)}\\
&\quad +y_{xx} y^{(4)}-y_x y^{(5)}+y y^{(6)}.
\end{split}
\end{align}
It is visible by inspection that these first integrals reduce for $\fq=0$ to those found in  \cite[(4.15)]{part1}  for trivial equations. An extension  of Corollary \ref{cor:gen-vsym} to the whole class of nonlinear equations of maximal symmetry can be formulated as follows,  and this extension also includes a result for first integrals.
\begin{cor}\label{cor:vsym-nlodes}
Let $\Delta_n[y]= 0$ be an $n$th order nonlinear {\small \sc ode} of maximal symmetry, and assume that it is obtained from the trivial equation by a point transformation $\sigma$ of the form \eqref{pt}, and denote by $\sigma_*$ the push-forward of $\sigma.$ Denote by  $\mcs_{div}$ the divergence symmetry algebra of the nonlinear {\small \sc ode}. In case this equation is also Lagrangian, denote by $\mcs_{var}$ its variational symmetry algebra.
\begin{enumerate}
\itemsep=1.5mm

\item[(a)] For $n$ odd, $\mcs_{div}$ is isomorphic under $\sigma_*$ to $\mathcal{A}_n \dotplus \F W_y.$
\item[(b)] For $n$ even, $\mcs_{div}$ is isomorphic under $\sigma_*$ to $\mathcal{A}_n \dotplus \mathfrak{s}.$
\item[(c)] If $\Delta_n[y]= 0$ is a Lagrangian equation, then with respect to the Lagrangian $\mscl_n$ given by \eqref{lagran-trf}, $\mcs_{var}$ is isomorphic under $\sigma_*$  to $\langle V_0, \dots,  V_{\frac{n-2}{2}}\rangle \dotplus  \langle F_n, G_n\rangle.$
\item[(d)] A fundamental set of $n$ linearly independent first integrals of the nonlinear equation $\Delta_n[y]=0$ is obtained from the first integral  \cite[(4.13)]{part1}  by a mere change of coordinates defined by $\sigma.$
\end{enumerate}
\end{cor}
\begin{proof}
Parts (a), (b) and (c) of the corollary are just a direct application of Corollary \ref{cor:gen-vsym} and Theorem \ref{th:transf-vsym}, and the well known fact due to Lie \cite{LieCanonic, olv94, olv-EIS} that any nonlinear {\small \sc ode} of maximal symmetry can be reduced to the canonical form by a point transformation. This result of Lie about the existence of the point transformation also implies part (d) of the corollary, and completes the proof.
\end{proof}

\begin{example}\label{ex:nleq}
The nonlinear {\small \sc ode}
\begin{equation}\label{eq:nleq}
\frac{1}{y^4} \left( 6 y_x^4- 12 y\,  y_x^2\, y_{xx}+  3 y^2  y_{xx}^2 + 4 y^2  y_x \,y^{(3)} -y^3 y^{(4)} \right)=0
\end{equation}
can be verified to be a Lagrangian equation of maximal Lie point symmetry algebra, of dimension $8.$ Therefore, the structures of its divergence and variational symmetry algebras are completely determined without any additional information by Corollary \ref{cor:vsym-nlodes}. In particular the corollary indicates that the divergence symmetry algebra has dimension $7$ while the variational symmetry algebra has dimension $4.$ The quickest way to find explicit expressions for the generators of these symmetry algebras as well as that for its first integrals is to find the point transformation $\sigma = \sigma(z, w)$ that maps the trivial equation to \eqref{eq:nleq}. It turns out that such a transformation $\sigma$ is given by
\begin{equation}\label{eq:trans-nl}
z=x,\qquad w= k_2 - k_1 \ln(y),
\end{equation}
where $k_1$ and $k_2$ are arbitrary constants.  Denoting by $X^t_j$ the transformed version under \eqref{eq:trans-nl} of a symmetry generator $X_j$ of the trivial equation, it follows from Corollary \ref{cor:vsym-nlodes} that for \eqref{eq:nleq} the divergence symmetry algebra $\mcs_{div}$ is given by
\[ \mcs_{div} = \langle V^t_0, V^t_1, V^t_2, V^t_3\rangle \dotplus  \langle F^t_4, G^t_4, H^t_4 \rangle
\]
where
\begin{alignat*}{3}
V^t_0 &= -y\pd_y,&\quad              V^t_1 &= -x y \pd_y,&\quad   V^t_2 &= -x^2 y \pd_y, \\
\quad  V^t_3 &= -x^3 y \pd_y,&\quad     G^t_4 &=2 x \pd_x -3y (k_2 - \ln(y)),&\quad         &  \\
F^t_4 &= \pd_x,&\quad        H^t_4 &= -x^2 \pd_x + 3 x y (k_2-\ln(y) ).&        &
\end{alignat*}
On the other hand, it also follows from Corollary \ref{cor:vsym-nlodes} that with respect to the transformed Lagrangian
$\mscl= -(y_x^2 - y y_{xx})^2/ (2 y^4)$, the variational symmetry algebra $\mcs_{var}$ of \eqref{eq:nleq} is given by
\[ \mcs_{var} = \langle V^t_0, V^t_1 \rangle \dotplus \langle F^t_4, G^t_4\rangle.\]

Finally, applying \eqref{eq:trans-nl} to  \cite[(4.13)]{part1}  shows that four linearly independent first integrals of \eqref{eq:nleq} are given by the linear combination
\[
 \begin{split}
\mathscr{F} = -\frac{1}{y^3} &\bigg[ a_0 \left(2 y_x^3-3 y y_x y_{xx}+y^2 y^{(3)}\right) + a_1 \bigg(-2 x y_x^3-y y_x \left(y_x-3 x y_{xx}\right)\\
&+y^2 (y_{xx}-x y^{(3)})\bigg) -a_3 \bigg(6 y^3 (k_2-\ln(y)) +2 x^3 y_x^3\\
&+3 x^2 y y_x \left(y_x-x y_{xx}\right)+x y^2 \left(6 y_x+x (-3 y_{xx}+x y^{(3)})\right)\bigg) \\
&+ a_2 \bigg(-2 x^2 y_x^3+x y y_x (-2 y_x+3 x y_{xx})\\
&-y^2 \left(2 y_x+x (-2 y_{xx}+x y^{(3)})\right)\bigg)
\bigg]
\end{split}
\]
where the arbitrary constants $a_j$ are  as in   \cite[(4.13)]{part1}  the coefficients of the linear combination of the individual first integrals $\mathscr{F}(V^t_j)$  associated with the symmetry generators $V^t_j$.
\end{example}

\section{On the variational symmetry algebra for the natural Lagrangian} \label{s:natlagrang}

 As already indicated, the most natural and most readily available Lagrangian for linear equations $\Delta_n[y]=0$ has a very simple closed form expression given by $\mcl_n= {\frac{1}{2}y \Delta_n[y]}.$ Recall that this $n$th order Lagrangian has been reduced to the Lagrangian $L_n$ of order $n/2$ of the form  \cite[(3.3)]{part1}  by making use of identities of the form  \cite[(3.2)]{part1} . However, although the variational symmetry algebra for \eqref{nor1} is the same relative to  $\mcl_n$ and $L_n,$ it is completely different relative to the transformed version $\mscl_n,$ as computed in full in the previous section.\par

Given the importance and the specific features of the natural Lagrangian $L_n,$ we shall also investigate the variational symmetry algebra of \eqref{nor1} relative to $L_n.$ However, in so doing we shall restrict our attention in this section to finding only those basis vectors of the divergence symmetry algebra  $\mathcal{A}_n \dotplus \mathfrak{s},$ as given by Corollary \ref{cor:gen-vsym}, which may be variational symmetries relative to the Lagrangian $L_n.$ We shall therefore be interested in this section only  in Lagrangian equations, and more precisely in equations of even orders  \cite[(3.1)]{part1} .

 We recall that the constant Wronskian of two linearly independent solutions $u, v$ of the source equation \eqref{srce} will be assumed to be normalized to 1. That is
\begin{equation} \label{wron=1}
\mathpzc{w} (u,v):= u v'-u'v=1.
\end{equation}
Moreover, assuming that $u$ and $v$ are linearly independent solutions of \eqref{srce} implies that
\begin{equation} \label{rul4a}
\fq = \mathpzc{w} (u',v'),\qquad \text{ but also}\qquad  \fq= -\frac{u''}{u}= -\frac{v''}{v}.
\end{equation}
In particular, it is always possible to replace derivatives of $u$ and $v$ of order two or higher by expressions involving derivatives of $u$ or $v$ of order at most one through the substitution
\begin{equation} \label{uvderiv}
D^{j+2} h = -D^j (\fq\, h),\qquad \text{ for $j \geq 0\;$ and $\;h=u$ or $v.$}
\end{equation}
It is also clear that one may always choose $u={\rm Const.}$ or $v={\rm Const.}$ if and only if $\fq=0.$
\subsection{The case $\bv=V_k,\quad 0 \leq k \leq n-1,\quad (n \in \set{4,6,8})$}

In order to investigate the vanishing of $\mathscr{S}(\mathbf{v}),$ we expand it as a differential polynomial in $y=y(x)$ and require that the resulting coefficients vanish identically. For $n=4,$ the vanishing of the coefficient of $y_x$ in this expansion corresponds to
\begin{equation} \label{n=4e1} 10 u^{2-k} v^{k-1} \fq\, b_0 =0,\qquad \text{ where } b_0= (k-3) v u_x - k u v_x,
\end{equation}
that is, to $\fq=0$ or $b_0=0.$ The condition $\fq=0$ is equivalent on account of \eqref{rul4a} to
\begin{equation} \label{u4a=0} u= \lambda + \theta v \end{equation}
for some scalars $\lambda$ and $\theta\neq 0.$ Moreover, it should be noted that the roles of $u$ and $v$ are symmetrical in the expression of the solution $s_k.$ Substituting \eqref{rul4a} and then \eqref{u4a=0} in the full equation $\mathscr{S}(\mathbf{v})=0$ representing \cite[(2.4)]{part1}  reduces this expression to
\begin{equation} \label{n=4e2}
v^{k-2} u^{1-k} \left[ (k-1) k \lambda^2 + 4 k \theta \lambda v + 6 \theta^2 v^2 \right]v_x^2 y_{xx}=0.
\end{equation}
Thus the linear independence of $1, v$ and $v^2$ shows that $k=0$ or $k=1$ and $\theta= 0,$ i.e. $u= \text{\rm Const.}$ Conversely for $k\in \set{0,1}$ and $u= \text{\rm Const.}$, we readily see that $\mathscr{S} (\mathbf{v})=0.$ Assuming that $\fq\neq 0,$  we have $\mathscr{S} (\mathbf{v})=0$ if and only if $b_0=0,$ and we readily see that there are no possible values for $k$ or $u$ for which the variational symmetry condition  \cite[(2.4)]{part1}  holds. In conclusion, when $n=4,$ $\mathbf{v}= V_k$  is a variational symmetry if and only if $k \in \set{0,1}$ and $u= \text{\rm Const.},$ or equivalently if and only if $v= \text{\rm Const.}$ and $k \in \set{2,3}.$ \par

A similar result holds for $n=6$ and $n=8.$ For these values of $n,$ the coefficients of $y_{xx}$  and $y^{(3)}$ in the expansion of $\mathscr{S}(\mathbf{v})$ are $(37 \fq^2 + 9 \fq_{xx})$ and $47 \fq \fq_x + 2 \fq^{(3)},$ respectively. The vanishing of these coefficients shows that all higher derivatives of $\fq$ can be expressed in terms of lower derivatives. The resulting expression for   \cite[(2.4)]{part1}  always shows that $\fq=0,$ and more precisely $u= \text{\rm Const.}$ and $k \in \set{0,1, \dots, (n-2)/2},$ and in turn these values of $u$ and $k$ lead to the vanishing of $\mathscr{S} (\mathbf{v}).$ We have thus established the following result.
\begin{thm}\label{th:vsym-Vk}
For $n \in \set{4,6,8},$ $\mathbf{v}= V_k$ is a variational symmetry for $\Delta_n[y]=0$ relative to $L_n$ if and only if
$\fq=0.$ In this case,  $k \in \set{0,1, \dots, (n-2)/2}$ if $u= {\rm Const.}$ or equivalently $k \in \set{n/2, \dots, n-1}$ if $v={\rm Const.}$
\end{thm}

\subsection{$\mathbf{v}$ is a generator of the $\mathfrak{sl}_2$-triplet $\frak{s}$}

We now consider the case where the symmetry generator $\mathbf{v}=F_n, G_n$ or $H_n,$ for $n \in \set{4,6,8}.$ We first let $n=4$ and $\mathbf{v}=F_4.$ Then as usual we expand $\mathscr{S}(\mathbf{v})$ as a polynomial in $y$ and its derivatives. The vanishing of the coefficients of $y_x y_{xx}$ in the resulting expression corresponds to the condition
\begin{equation} \label{F4e1}
\fq = u_x^2/ u^2,
\end{equation}
and the latter condition together with \eqref{rul4a} yields the system of equations
\begin{equation} \label{vsymcdt4F4}
u u_{xx} + u_{x}^2 =0, \qquad u v_x - u_x v=1.
\end{equation}
In the sequel, the $k_j, \lambda$ and $\alpha$ will denote some arbitrary constants.  Solving \eqref{vsymcdt4F4} yields, for nonconstant $u$ , the solutions
\begin{equation} \label{soln4F4}
u= k_2 \sqrt{2x - k_1},\qquad v= \frac{\sqrt{2x - k_1}}{2 k_2}\left[ 2 k_2^2 k_3 + \ln(k_1- 2x)  \right].
\end{equation}
Conversely it is easily verified that $F_4$ is a variational symmetry when $u$ is a constant or when $u$ and $v$ are given by \eqref{soln4F4}. Since the expression for $H_4$ is up to sign the same as that for $F4$ in which $u$ is replaced by $v,$ a similar reasoning shows that $H_4$ ia a variation symmetry if and only if $v$ is the constant function or if
\begin{equation} \label{soln4H4}
v= k_2 \sqrt{2x - k_1},\qquad u= \frac{\sqrt{2x - k_1}}{2 k_2}\left[ 2 k_2^2 k_3 - \ln(k_1- 2x)  \right].
\end{equation}

It should also be noted that for the expressions of $u$ and $v$ given in \eqref{soln4F4} and \eqref{soln4H4}, we have correspondingly
\begin{equation*}
\fq= \frac{1}{(k_1 -2 x)^2}.
\end{equation*}


We now move on to the case where $n=4$ and $\mathbf{v} =G_4.$ Here the vanishing of the coefficient of $y_x y_{xx}$ in the polynomial expansion of $\mathscr{S}(\mathbf{v})$ shows that
\begin{equation} \label{AforG4}
\fq = \frac{u_x v_x}{u v}.
\end{equation}
Using in particular the fact that $u$ and $v$ are linearly independent solutions of \eqref{srce} leads to the following system of equations
\begin{equation} \label{vsymcdt4G4}
\frac{u_x v_x}{u v}= - u_{xx}/u,\qquad u v_x-u_xv=1,
\end{equation}
which is reducible to the following one.
\begin{equation} \label{vsymcdt4G4v2}
v= \lambda/u_x, (u_x \neq 0), \qquad  u u_xx + \alpha u_x^2 =0,\qquad  \alpha = 1+ \frac{1}{\lambda}.
\end{equation}
Solving \eqref{vsymcdt4G4v2} gives rise to the following solution depending on the arbitrary parameter $\alpha.$
\begin{align} \label{soln4G4}
\begin{split}
u &= {\rm Const.} \text{ or } v=  {\rm Const.} \text{ (in particular if $\alpha =0$), or} \\
u &= \frac{k_2}{\sqrt{-\alpha}} e^{k_1 x}, \qquad v = \frac{\lambda \sqrt{-\alpha}}{k_1 k_2 e^{k_1 x}},\qquad  \text{ (if $\alpha <0$), or}\\
u & = \frac{k_2}{\sqrt{\alpha}} \sqrt{2 x - k_1},\qquad v = \frac{\lambda \sqrt{\alpha}}{k_2} (2 x- k_1)^{-3/2}, \quad \text{ (if $\alpha >0$)}.
\end{split}
\end{align}

In turn, it is easily verified that for the values of $u$ and $v$ in \eqref{soln4G4}, $\bv= G_4$ is a variational symmetry. For $n=6$ and $n=8$ a similar analysis yields exactly the same results obtained for the case $n=4.$
\begin{thm} \label{th:vsym-sl2}
Let $\Delta_n[y]=0$ be a {\small \sc lode} of maximal symmetry as in  \cite[(3.1)]{part1}, considered as a Lagrangian equation relative to $L_n,$ with $n \in \set{4,6,8}.$
\begin{enumerate}
\item[(a)] $\mathbf{v}= F_n$ is a variational symmetry if and only if either $u= {\rm Const.}$ or  else $u$ and $v$ are given by \eqref{soln4F4}.
\item[(b)] $\mathbf{v} = G_n$ is a variational symmetry if and only if either $\fq=0$ or else $u$ and $v$ are given by \eqref{soln4G4}.
\item[(c)] $\mathbf{v}= H_n$ is a variational symmetry if and only if either $v= {\rm Const.}$ or else $u$ and $v$ are given by \eqref{soln4H4}.
\end{enumerate}
\end{thm}
Recall from Corollary \ref{cor:gen-vsym} that relative to the transformed Lagrangian $\mscl_n,$ the whole variational symmetry algebra of \eqref{nor1} was found to be $\langle V_0, \dots,  V_{\frac{n-2}{2}}\rangle \dotplus  \langle F_n, G_n\rangle.$ It therefore follows from Theorem \ref{th:vsym-Vk} and Theorem \ref{th:vsym-sl2} that not only the variational symmetry algebra of \eqref{nor1}  relative to $\mscl_n$ on one hand,  and to $L_n$ on the other hand are completely different, but also that this algebra is severely reduced when considered relative to $L_n.$ In particular, while for any value of $\fq$ half of the solution symmetries generate the variational symmetry algebra relative to $\mscl_n$, the variational symmetry algebra relative to $L_n$ is nontrivial if and only if $\fq=0.$\par

More generally, variational symmetry algebras of \eqref{nor1} relative to $\mscl_n$ which are valid for arbitrary values of $\fq$ will become reduced in size relative to $L_n$ and moreover, such reduced algebras will also be valid only for severely restricted values of $\fq$ which may  depend only on one arbitrary constant.\par

This independence concerning of variational symmetry algebras of equivalent Lagrangian equations which does not seem to have been discussed in the literature can however be easily explained. Let $L_0$ and $L$ be two Lagrangians corresponding to two equivalent Lagrangian equations. Thus $L= \theta L_0 + \divgce P$ for a certain differential function $P$ and a constant $\theta.$ For any vector field $\bv$ in the domain of the operator $\mscs$ defined in \eqref{plv}, let us write the expression $\mscs (\bv)$  as $\mscs_L(\bv)$ to underline the fact that the lagrangian being considered is $L.$  It then follows that
\begin{equation} \label{S_l(v)}
\mscs_L(\bv)= \theta \mscs_{L_0} (\bv) + \mscs_{\divgce P} (\bv).
\end{equation}
Now, given $L_0,$ there is an infinity of arbitrary choices for $\divgce P,$ each of which gives rise to an expression for $L$ up to an arbitrary constant.  Therefore, \eqref{S_l(v)} clearly shows that a variational symmetry $\bv$ for the Lagrangian equation corresponding to $L_0$ needs not also be a variational symmetry for the equivalent lagrangian equation corresponding to $L,$ because the expression $\mscs_{\divgce P} (\bv)$ will not vanish in general, even when $\mscs_{L_0} (\bv)$ does. This is precisely what happens with $L= L_n$ and $L_0 = \mscl_n,$  as given by  \cite[(3.3)]{part1}  and \eqref{lagran-trf},  respectively.

\section{Concluding Remarks} \label{s:conclu}
One of the main results of this paper is Corollary \ref{cor:vsym-nlodes} which completely determines the structure of the variational and the divergence symmetry algebras of any nonlinear {\small \sc ode} of maximal symmetry based solely on its order. This corollary also determines a fundamental set of first integrals of the nonlinear equation once the linearizing transformations mapping it to the canonical form are known. All these results obtained for nonlinear {\small \sc ode}s are generally based on those obtained in \cite{part1} for linear equations in canonical form. \par

For two systems of equivalent linear or nonlinear Lagrangian equations, the relationship established in \eqref{S_l(v)} shows that not only equivalent Lagrangian equations as defined in the paper need not have isomorphic variational symmetry algebras, but also that the variational symmetry algebra of one of the Lagrangian equations can be chosen independently of the other. However, the complete description of  this independence  appears to be an open problem. Section \ref{s:natlagrang} was devoted to an example of such an independence between the variational symmetry algebras of equivalent Lagrangian equations associated with the Lagrangians $L_n$ and $\mscl_n,$ respectively.

\end{document}